\documentclass[12pt,reqno]{amsart}
\usepackage{amscd,amssymb,amsthm}
\usepackage{graphicx}
\usepackage{enumerate}
\theoremstyle{plain}
\newtheorem{theorem}{Theorem}[section]
\newtheorem{corollary}[theorem]{Corollary}

\newcommand{\abs}[1]{\left\vert #1\right\vert }

\newcommand{\bg}{\medskip\goodbreak}
\newcommand{\itemref}[1]{\eqref{#1}}
\newenvironment{enumeratea}{\begin{enumerate}%
	[\upshape (a)]}{\end{enumerate}}
\newenvironment{enumeratei}{\begin{enumerate}%
	[\itshape i.]}{\end{enumerate}}
\title[The Positivity of the Coefficients of Some Series Expansions]
{A Note on The Positivity of the Coefficients of Some\\
 Power Series Expansions}
\author[Omran Kouba]{Omran Kouba$^\dag$}
\address{Department of Mathematics \\
Higher Institute for Applied Sciences and Technology\\
P.O. Box 31983, Damascus, Syria.}
\email{omran\_kouba@hiast.edu.sy}
\thanks{$^\dag$ Department of Mathematics, Higher Institute for Applied Sciences and Technology.}

\begin{document}
\parindent=0pt
\date{\today}
\begin{abstract}
In this short note, a general result concerning the positivity, under some conditions, of the coefficients of a power
series is proved. This allows us to answer positively a question raised by Guo (2010)  about
the sign of the coefficients of a power series relating the residual errors in Halley's iterations for the $p$th root.
\end{abstract}
\smallskip\goodbreak

\maketitle

\section{\bf Introduction }\label{sec1}
\bg
\parindent=0pt
\qquad
The determination of the sign of the coefficients of power series expansions of rational functions
was considered several times in the literature (see \cite{kra},\cite{str} and the bibliography therein).\par
\qquad Let us describe the context in which our result is situated. Let $p\geq 2$ be an integer, and let  $z$ be any complex number. If we apply Newton's method to solve the equation $x^p=1-z$ starting from the initial value $1$, we get
the sequence of rational functions $(U_k)_{k\geq0}$ in the variable $z$ defined by the iteration

\begin{equation}\label{E:eqn}
U_{k+1}(z)=\frac{1}{p}\left((p-1)U_k(z)+\frac{1-z}{U_k^{p-1}(z)}\right),\qquad U_0(z)\equiv 1.
\end{equation}

\qquad Similarly, if we apply Halley's method to solve the equation $x^p=1-z$ starting from the same initial value $1$, we get
the sequence of rational functions $(V_k)_{k\geq0}$ in the variable $z$ defined by the iteration

\begin{equation}\label{E:eqh}
V_{k+1}(z)=\frac{(p-1)V_k^p(z)+(p+1)(1-z)}{(p+1)V_k^{p}(z)+(p-1)(1-z)}V_k(z),\qquad V_0(z)\equiv 1
\end{equation}
\bg
\qquad Following \cite{bini} and \cite{guo}, we define, for Newton's method or Halley's method,  the residual errors : $(N_k)_{k\geq0}$ and $(H_k)_{k\geq0}$
by
$$N_k(z)=1-\frac{1-z}{U_k^p(z)} ,\qquad\hbox{and}\qquad H_k(z)=1-\frac{1-z}{V_k^p(z)},$$

and one checks easily that
\begin{equation}\label{E:eqfg}
N_{k+1}(z)=f_p(N_k(z))\qquad\hbox{and}\qquad H_{k+1}(z)=g_p(H_k(z)), 
\end{equation}
where,
\begin{align}
f_p(t)&=1-(1-t)\left(1-\frac{t}{p}\right)^{-p}, \label{E:eq00}\\
g_p(t)&=1-(1-t)\left(\frac{2p-(p-1)t}{2p-(p+1)t}\right)^{p}. \label{E:eq0}
\end{align}
\qquad In \cite{guo} it was shown that for $w\in\overline{D(0,1)}\setminus\{0,1\}$ we have
$\abs{f_p(w)}<\abs{w}^2$.
 Hence,
the conditions $\abs{N_k(z)}\leq1$ and $N_k(z)\notin\{0, 1\}$ imply that $\abs{N_{k+1}(z)}<\abs{N_k(z)}^2$, 
and consequently, for every $k\geq 1$ we have
\bg
\begin{equation}\label{E:eq4}
\abs{N_k(z)}< \abs{z}^{2^k},\qquad\hbox{for $\abs{z}\leq 1,~ z\notin\{0, 1\}$ .}
\end{equation}

This is essentially Lemma 1 from \cite{guo}.  In our Theorem \ref{th23}, we will prove a similar result about the sequence $(H_k(z))_{k\geq0}$ of 
residual errors corresponding to Halley's method. This will answer positively a question asked in \cite{guo}.
\bg

\section{\bf The Main Results }\label{sec2}
\bg
\qquad We start this section by proving a general result :
\bg
\begin{theorem}\label{th21}
Consider a non-constant monotone decreasing sequence of positive real numbers   $(a_n)_{n\geq1}$ with $a_1=1$, and define $\ell =\min\{k\geq 1:a_k<1\}$. Let $F$ be the function defined in the open unit disk $D(0,1)$ by the formula :
$$F(z)=1-(1-z)\exp\left(\sum_{n=1}^\infty\frac{a_n}{n}z^n\right),$$
and consider  the power series expansion $\sum_{n=1}^\infty c_nz^n$ of $F$ in the neighbourhood of $0$. Then
\begin{enumeratea}
       \item  For every $n\geq0$ we have $c_n\in[0,1)$, and $\sum_{n=1}^\infty c_n\leq 1$. \label{itm1}
       \item  $F$ can be  continuously extended to the closed unit disk $\overline{D(0,1)}$. \label{itm2}
       \item  For $n<\ell$ we have $c_n=0$, and  for $n\geq \ell$, we have $c_n>0$.\label{itm3}
       \item  For $z\in\overline{D(0,1)}\setminus\{0,1\}$ we have $\abs{F(z)}<\abs{z}^\ell$.\label{itm4}
\end{enumeratea}
\end{theorem}
\bg
\begin{proof}
Indeed, let us define $G$ and $H$ in $D(0,1)$ by
\begin{equation}
          G(z)=\sum_{n=1}^\infty\frac{a_n}{n}z^n,\qquad \hbox{and}\qquad H(z)=e^{G(z)},
\end{equation}
and let us suppose that $\sum_{n=0}^\infty b_nz^n$ is the power series expansion of $H$. Clearly  $b_0=1$, and
\begin{equation*}
         (1-z)H(z)=1+\sum_{n=0}^\infty b_{n+1}z^{n+1}-\sum_{n=0}^\infty b_{n}z^{n+1}
\end{equation*}
\bg
so,
\begin{equation*}
              F(z) =1-(1-z)H(z)=\sum_{n=0}^\infty (b_{n}-b_{n+1})z^{n+1},
\end{equation*}
and consequently 
\begin{equation}\label{E:eq210}
               \forall\, n\geq0,\qquad c_{n+1}=b_{n}-b_{n+1}
\end{equation}

On the other hand, noting that $H'(z)=G'(z)H(z)$ we conclude that
\begin{equation*}
\sum_{n=0}^\infty (n+1)b_{n+1}z^n=\left(\sum_{n=0}^\infty a_{n+1} z^{n}\right)\left(\sum_{n=0}^\infty b_nz^n\right).
\end{equation*}
This shows that for $n\geq 0$ we have
\begin{equation}\label{E:eq211}
b_{n+1}=\frac{1}{n+1}\sum_{k=0}^{n}a_{k+1}b_{n-k}.
\end{equation}
In particular, $b_1=1$ and  $b_n>0$ for every $n\geq0$.\bg
Recalling that $a_1=1$ we conclude from \eqref{E:eq211} that
\begin{equation}\label{E:eq212}
       (n+1)b_{n+1}-b_n=\sum_{k=1}^{n}a_{k+1}b_{n-k},
\end{equation}
and 
\begin{equation}\label{E:eq213}
      n b_{n}=\sum_{k=0}^{n-1}a_{k+1}b_{n-1-k}=\sum_{k=1}^{n}a_{k}b_{n-k}.
\end{equation}
So, subtracting  \eqref{E:eq212} from \eqref{E:eq213} and rearranging, we obtain 

\begin{equation}\label{E:eq214}
     b_n-b_{n+1}=\frac{1}{n+1}\sum_{k=1}^{n}(a_k-a_{k+1})b_{n-k},
\end{equation}

which is valid for $n\geq1$. (It would be valid also for $n=0$ if  we interpret the right hand side sum as 0 in that case.)  \bg
We deduce from \eqref{E:eq210} and \eqref{E:eq214} that $c_1=0$ and 

\begin{equation}\label{E:eq215}
        \forall\, n\geq2,\qquad c_{n}=\frac{1}{n}\sum_{k=2}^{n}(a_{k-1}-a_{k})b_{n-k},
\end{equation}

Now, since the sequence $(a_k)_{k\geq1}$ is decreasing, and  the $b_m$'s are positive, we conclude that $c_n\geq 0$ for every $n\geq1$.  It follows also from \eqref{E:eq210} that 

\begin{equation*}
        \forall\,m\geq1,\qquad\sum_{n=1}^m c_n=1-b_m,
\end{equation*}
this proves, in particular, that $c_n\in[0,1)$ for all $n\geq1$, and that the series $\sum_{n\geq1}c_n$ is convergent with sum smaller or equal to $1$. Hence, the power series expansion of $F$ is normally convergent in the closed unit disk $\overline{D(0,1)}$. This 
concludes the proof of  both \itemref{itm1} and \itemref{itm2}.\bg
\qquad To see \itemref{itm3} we note that
$$1=a_1=a_2=\ldots=a_{\ell-1}>a_\ell,$$
therefore by \eqref{E:eq215}  if $n <\ell$ then $c_n=0$, and if $n\geq\ell$ then
$c_n\geq (a_{\ell-1}-a_\ell)b_{n-\ell}/n>0$. \bg
\qquad Finally, for $z\in\overline{D(0,1)}\setminus\{0,1\}$ we have
\begin{align*}
\frac{\abs{F(z)}}{\abs{z}^\ell}&=\abs{\sum_{n=\ell}^\infty c_nz^{n-\ell}}\leq\abs{c_\ell+z c_{\ell+1}} +\sum_{n=\ell+2}^\infty c_n\notag\\
&<c_\ell+ c_{\ell+1} +\sum_{n=\ell+2}^\infty c_n\leq 1.
\end{align*}
This completes the proof of \itemref{itm4} and achieves the proof of the theorem.
\end{proof}
\bg
{\bf Examples :}\par
\begin{enumeratei}
\item  Let $p>1$ be a real number, and consider, for $n\geq1$, $a_n=p^{1-n}$ . Applying Theorem \ref{th21} to this
data proves that the function $f_p$ defined on $D(0,p)$ by
$$f_p(z)=1-(1-z)\left(1-\frac{z}{p}\right)^{-p},$$
satisfies the inequality $\abs{f_p(z)}<\abs{z}^2$ for $z\in\overline{D(0,1)}\setminus\{0,1\}$.\label{ex1}
\item  Let $p>1$ be a real number, and consider, for $n\geq1$, $a_n=p(\alpha^n-\beta^n)$
where $\alpha=\frac{p+1}{2p}$ and $\beta=\frac{p-1}{2p}$. Clearly, $a_n\geq0$ because $\alpha>\beta>0$, and
$$a_1=a_2=1,\qquad  a_n-a_{n+1}=\alpha\beta a_{n-1}>0\quad\hbox{for $n\geq2$}.$$
Applying Theorem \ref{th21} to this
data proves that the function $g_p$ defined on $D(0,\frac{2p}{p+1})$ by
$$g_p(z)=1-(1-z)\left(\frac{2p-(p-1)z}{2p-(p+1)z}\right)^p,$$
satisfies the inequality $\abs{g_p(z)}<\abs{z}^3$ for $z\in\overline{D(0,1)}\setminus\{0,1\}$.\label{ex2}
\end{enumeratei}
\bg
\qquad The conclusion of Example {\itshape\ref{ex1}}. was used in \cite{guo} to prove the following corollary about the residual errors in Newton's method, (see Lemma 1 in \cite{guo}) :\bg
\begin{corollary}\label{cor22}
Let $p\geq 2$ be an integer, and let $(U_k)_{k\geq0}$ be the sequence of complex rational functions defined by 
the iteration
$$
U_{k+1}(z)=\frac{1}{p}\left((p-1)U_k(z)+\frac{1-z}{U_k^{p-1}(z)}\right),\qquad U_0(z)\equiv 1.
$$
$($which is obtained when Newton's method is applied to solve $x^p=1-z$ starting from $1.)$ Then for every
$k\geq1$ and every $z\in \overline{D(0,1)}\setminus\{0,1\}$ we have
$$
\abs{1-\frac{1-z}{U_k^p(z)}}<\abs{z}^{2^k}.
$$
In particular, $U_k$ has no zeros, nor poles in $\overline{D(0,1)}$.
\end{corollary} 
\bg
\qquad In the same spirit, using the conclusion of Example {\itshape\ref{ex2}}.  we will prove a similar result about the residual errors in Halley's method,
which illustrates the cubic character of the convergence of Halley's method. It should be compared with Lemma 2 from \cite{guo}. 
\bg
\begin{theorem}\label{th23}
Let $p\geq 2$ be an integer, and let $(V_k)_{k\geq0}$ be the sequence of complex rational functions defined by 
the iteration
$$
V_{k+1}(z)=\frac{(p-1)V_k^p(z)+(p+1)(1-z)}{(p+1)V_k^{p}(z)+(p-1)(1-z)}V_k(z),\qquad V_0(z)\equiv 1,
$$
$($which is obtained when Halley's method is applied to solve $x^p=1-z$ starting from $1.)$ Then for every
$k\geq1$ and every $z\in \overline{D(0,1)}\setminus\{0,1\}$ we have
$$
\abs{1-\frac{1-z}{V_k^p(z)}}<\abs{z}^{3^k}.
$$
In particular, $V_k$ has no zeros, nor poles in $\overline{D(0,1)}$.
\end{theorem} 
\bg
\begin{proof}
Recalling \eqref{E:eqfg} and using the conclusion of Example {\itshape\ref{ex2}}. we see that if Halley's residual $H_k(z)$ satisfies the
the conditions $\abs{H_k(z)}\leq1$ and $H_k(z)\notin\{0, 1\}$, then $\abs{H_{k+1}(z)}<\abs{H_k(z)}^3$. Consequently, for every $k\geq 1$ we have
\begin{equation*}
\abs{H_k(z)}< \abs{z}^{3^k},\qquad\hbox{for $\abs{z}\leq 1,~ z\notin\{0, 1\}$ .}
\end{equation*}
Clearly, this implies that $V_k$ has no zeros, nor poles in $\overline{D(0,1)}$.
\end{proof}
\bg
\qquad It is easy to show by induction that $U_k(0)=1$ and $V_k(0)=1$ for every $k$. So, if $z\mapsto \root{p}\of{1-z}$ is the principal determination of the $p$th root, then we conclude from the following identities :
\begin{align*}
U_k(z)-\root{p}\of{1-z}=\frac{U_k^p(z)}{\sum_{r=1}^{p-1}U_k^{r-1}(z)\left(\root{p}\of{1-z}\right)^{p-r}}\,N_k(z),\\
\noalign{\hbox{and}}\\
V_k(z)-\root{p}\of{1-z}=\frac{V_k^p(z)}{\sum_{r=1}^{p-1}V_k^{r-1}(z)\left(\root{p}\of{1-z}\right)^{p-r}}\,H_k(z),
\end{align*}
that, in the neighbourhood of $z=0$, we have
\begin{equation}\label{E:eq200}
U_k(z)-\root{p}\of{1-z}=O(z^{2^k}),
\quad\hbox{and}\quad
V_k(z)-\root{p}\of{1-z}=O(z^{3^k}),
\end{equation}
for every $k\geq0$. 
\bg
\qquad In fact, \eqref{E:eq200} provides a different proof of Theorem 10 in \cite{guo}, recalled in a different form in the
following Corollary.
\bg
\begin{corollary}\label{cor23}
Let $p\geq 2$ be an integer, and let $(U_k)_{k\geq0}$, and $(V_k)_{k\geq0}$ be the sequences of complex rational functions defined by the iterations \eqref{E:eqn} and \eqref{E:eqh}. Then, for every $k\geq0$, the
rational functions $U_k$ and $V_k$
have power series expansions converging in a neighbourhood of the closed unit disk, and
\begin{align*}
U_{k}(z)=1-\sum_{n=1}^{2^k-1}\left(\prod_{r=1}^{n-1}(rp-1)\right)\frac{z^n}{n!\,p^n}+O(z^{2^k}),\\
\noalign{\hbox{and}}\\
V_{k}(z)=1-\sum_{n=1}^{3^k-1}\left(\prod_{r=1}^{n-1}(rp-1)\right)\frac{z^n}{n!\,p^n}+O(z^{3^k}),
\end{align*}
where $\prod_{r=1}^{n-1}(rp-1)$ is interpreted as $1$ for $n=1$.
\end{corollary} 
\bg


\end{document}